\documentclass[11pt]{article}
\usepackage{mathrsfs}
\usepackage{amsmath}
\usepackage{amssymb}
\usepackage{graphicx}
\usepackage{epic}
\usepackage{float}

\renewcommand{\paragraph}{\roman{paragraph}}

\def \la{\lambda}

\newtheorem{theorem}{\scshape \mdseries  Theorem}[section]
\newtheorem{lemma}[theorem]{\scshape \mdseries  Lemma}
\newtheorem{coro}[theorem]{\scshape \mdseries  Corollary}

\begin{document}

\title{\sf  Lower bounds for algebraic connectivity of graphs in terms of matching number or edge covering number\thanks{
Supported by National Natural Science Foundation of China (11071002, 11371028, 71101002), Program for New Century Excellent
Talents in University (NCET-10-0001), Key Project of Chinese Ministry of Education (210091),
Specialized Research Fund for the Doctoral Program of Higher Education (20103401110002),
Natural Science Research Foundation of Department of Education of Anhui Province(KJ2012B040, KJ2013A196),
Scientific Research Fund for Fostering Distinguished Young Scholars of Anhui University(KJJQ1001).}
}
\author{Jing Xu$^1$,  Yi-Zheng Fan$^{1,}$\thanks{Corresponding author.
 E-mail addresses: fanyz@ahu.edu.cn(Y.-Z. Fan), xujing@ahu.edu.cn (J. Xu), tansusan1@aiai.edu.cn (Y.-Y. Tan).}, \ Ying-Ying Tan$^2$\\
    {\small  \it $1$. School of Mathematical Sciences, Anhui University, Hefei 230601, P. R. China}\\
    {\small  \it $2$. Department of Mathematics and Physics, Anhui University of Architecture, Hefei 230601, P. R. China} \\
  }
\date{}
\maketitle

\noindent {\bf Abstract:}
In this paper we characterize the unique graph whose algebraic connectivity is minimum among all connected
graphs with given order and fixed matching number or edge covering number,
and present two lower bounds for the algebraic connectivity in terms of the matching number or edge covering number.



\section{Introduction}

Let $G$ be a connected simple graph of order $n$ with vertex set
$V=V(G)=\{v_1,v_2,\ldots,v_n\}$ and edge set $E=E(G)$. The {\it
adjacency matrix} of the graph $G$ is defined to be the matrix
$A=A(G)=[a_{ij}]$, where $a_{ij}=1$ if $v_i$ is adjacent
to $v_j$, and $a_{ij}=0$ otherwise.
The matrix $L(G)=D(G)-A(G)$ is called the {\it Laplacian matrix} of the graph $G$, where
$D(G)=\hbox{diag}\{d_G(v_1), d_G(v_2),\cdots,d_G(v_n)\}$ is a diagonal
matrix, and $d_G(v)$ denotes the degree of a vertex $v$ in the graph $G$. It is
easy to see that $L(G)$ is a real and positive semidefinite, so that its
eigenvalues can be arranged as follows:
$$0 = \mu_n(G) \le \mu_{n-1}(G) \le \cdots \le \mu_1(G),$$
where $\mu_n(G)=0$ as each row sum of $L$ is zero, with the all-one vector $\mathbf{1}$ as an corresponding eigenvector. It is well known
that the multiplicity of eigenvalue 0 is equal to the number of
components of $G$. The eigenvalue $\mu_{n-1}(G)$, also denoted by $\alpha(G)$, is called
the {\it algebraic connectivity} of $G$ by Fiedler \cite{fied1}; and the eigenvectors
corresponding to $\alpha(G)$ are usually called the {\it Fiedler vectors} of $G$.

 The algebraic connectivity has received much attention; see \cite{abr,bappm,fallatk,Grone1,Grone2,kirk,merris,merris1,merris2,merris3,merris4}.
For example, upper bounding or  maximizing the algebraic connectivity has been discussed by
Lu et. al \cite{lu} in terms of the domination number, Lal et. al \cite{lal} subject to the number of pendant vertices,
Zhu \cite{zhu} by means of matching number.
Lower bounding or minimizing the algebraic connectivity has also been discussed by
Fallat et. al \cite{fallatk, fallatkp1} subject to diameter or girth,
Biyiko\v{g}lu and Leydold \cite{biy, biy2} subject to degree sequence or size,
and Fan and Tan \cite{fant} subject to domination number.

Recently Fan and Tan \cite{fant} obtain a perturbation result for the algebraic connectivity of a graph
 when a branch of the graph is relocated from one vertex to another vertex.
The result motivates us to do a lot of work on minimizing the algebraic connectivity subject to graph parameters, which
provides some lower bounds for the  algebraic connectivity.
In this paper, we characterize the unique graph whose algebraic connectivity is minimum among all connected
graphs with given order and fixed matching number or edge covering number,
and present two lower bounds for the algebraic connectivity in terms of the matching number or edge covering number.

At the end of this section, we introduce some notions.
Recall that a {\it matching} of a graph $G$ is an set of independent edges of $G$;
and the {\it matching  number} of $G$ is the maximal cardinalities of all the matchings of $G$, denoted by  $\beta(G)$.
Clearly, $n\geq 2\beta(G)$. In particular, $G$ has {\it perfect matchings}
 if $n=2\beta(G)$.
 An {\it edge cover} of a graph $G$ without isolated vertices is a set of edges of $G$ that covers all vertices of
$G$.
The {\it edge covering number} of a graph $G$ is the minimum cardinality
of all edge covers of $G$, denoted by $\gamma(G)$.
It is known that $\beta(G)+\gamma(G)=|V(G)|$ if $G$ contains no isolated vertices \cite{Gallai}.

Denote by $\mathcal{M}_{n,\beta}$ (respectively, $\mathcal{C}_{n,\gamma}$) the set of connected graphs of order $n$ with matching number $\beta$ (respectively, edge covering number $\gamma$).
Let $P_{d}$ denote a path of order $d$, and $S_{1,m}$ a star on $m+1$ vertices.
Denote by $T(k,l,d)$ a tree obtained from a path $P_{d}$
by attaching two stars $S_{1,k},S_{1,l}$ at its two end points respectively; see Fig. 1.1.
In particular, if $d=1$, then $T(k,l,d)=S_{1,k+l}$; if $k=1$ and $l=0$, then $T(k,l,d):=P_{d+1}$.
For convenience, a graph is called {\it minimizing} in a certain class if $\alpha(G)$ is minimum among all graphs in the class.%\in\mathcal{M}_{n,\beta}$.

\begin{center}
\includegraphics[scale=.6]{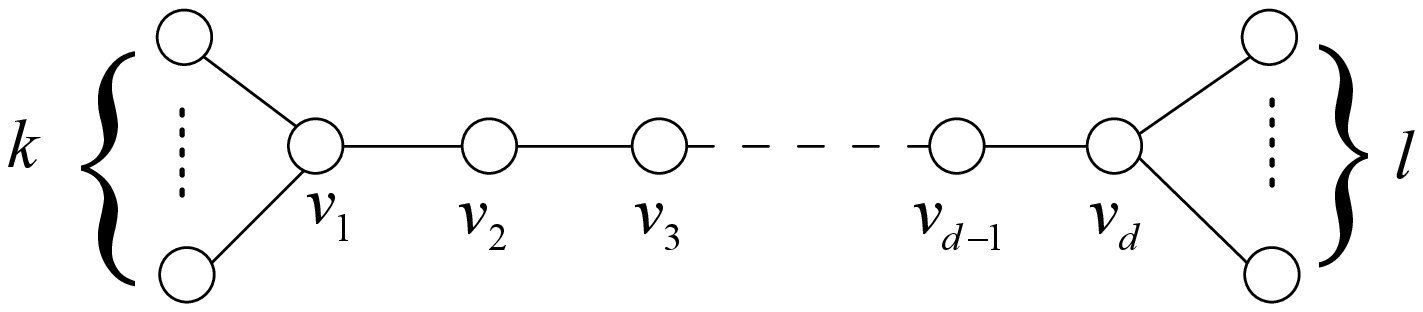}

{\small Fig. 1.1  The tree $T(k,l,d)$}
\end{center}

%Let $G$ be a graph on vertices $v_1,v_2,\ldots, v_n$, and
Let $x=(x_1, x_2, \ldots, x_n) \in \mathbb{R}^n$ and let $G$ be a graph on vertices $v_1,v_2,\ldots,v_n$.
The vector $x$ can be considered
as a function defined on $V(G)$, which maps each vertex $v_i$ of $G$ to the value $x_i$, i.e. $x(v_i)=x_i$.
If $x$ is an eigenvector of $L(G)$, then it defines on $G$ naturally, i.e. $x(v)$ is the entry of $x$ corresponding to $v$.
One can find
that the quadratic form $x^TL(G)x$ can be written as
$$x^TL(G)x=\sum_{uv \in E(G)}[x(u)-x(v)]^2. \eqno(1.1)$$
The eigenvector equation $L(G)x=\la x$ can be interpreted as
$$[d_G(v)-\la] x(v)= \sum_{u \in N_G(v)} x(u), \hbox{~ for each~} v \in V(G), \eqno(1.2)$$
where $N_G(v)$ denotes the neighborhood of $v$ in $G$.
In addition, for an arbitrary unit vector $x \in \mathbb{R}^n$ orthogonal to $\mathbf{1}$,
$$\alpha(G)\leq x^TL(G)x, \eqno(1.3)$$
with equality if and only if $x$ is a Fiedler vector of $G$.

%\begin{figure}[H]
%  \centering

% \includegraphics[width=4.5in]{fig1.eps}
%  \caption{\small The tree $T(k,l,d)$}\label{tree}
%\end{figure}

%If $d(v)=1$ for $v\in V(G)$, then $v$ is said to be {\it pendent vertex.}
%A vertex $v\in V(G)$ is called {\it quasi-pendant} if it is incident with
%a pendant vertex. Denote $K_{1,n-1}$ be the star graph with the root vertex
%$v$ adjacent to all the other vertices. If a star is an induced subgraph of
% some graph $G$, then we call it {\it pendent star} of graph $G$.

\section{Preliminaries}
First we introduce the property of Fiedler vectors of a tree.

\begin{lemma}{\em \cite{fied2}} \label{fiedler2}
Let $T$ be a tree with a Fiedler vector $x$.
Then exactly one of the two cases occurs:

Case A. All values of $x$ are nonzero. Then $T$ contains exactly one edge $pq$ such that $x(p)>0$ and $x(q)<0$.
The values in vertices along any path in $T$ which starts in $p$ and does not contain $q$ strictly increase, the values in vertices
along any path starting in $q$ and not containing $p$ strictly decrease.

Case B. The set $N_0=\{v: x(v)=0\}$ is non-empty. Then the graph induced by $N_0$ is connected and there is exactly one vertex $z \in N_0$ having
at least one neighbor not belonging to $N_0$. The values along any path in $T$ starting in $z$ are strictly increasing, or strictly decreasing, or zero.
\end{lemma}

If the Case B in Lemma \ref{fiedler2} occurs, the vertex $z$ is called the {\it characteristic vertex}, and $T$ is called a {\it Type I} tree;
otherwise, $T$ is called a {\it Type II} tree in which case the edge $pq$ is called the {\it characteristic edge}.
The characteristic vertex or characteristic edge of a tree is independent of the choice of Fiedler vectors; see \cite{merris}.

Next we introduce the perturbation result of the algebraic connectivity of a graph.
Let $G_1$, $G_2$ be two vertex-disjoint graphs, and let $v\in V(G_1)$, $u\in V(G_2)$. The {\it coalescence} of
$G_1$ and $G_2$ with respect to $v$ and $u$, denoted by $G_1(v)\diamond G_2(u)$, is obtained from $G_1$ and
$G_2$ by identifying $v$ with $u$ and forming a new vertex $p$, which is also denoted as $G_1(p)\diamond G_2(p)$.
If a connected graph $G$ can be expressed as $G=G_1(p)\diamond G_2(p)$,
where $G_1$ and $G_2$ are nontrivial subgraphs of $G$ both containing $p$,
then $G_1$ or $G_2$ is called a {\it branch} of $G$ rooted at $p$.
Let $G=G_1(v_2)\diamond G_2(u)$ and $G^*=G_1(v_1)\diamond G_2(u)$, where $v_1$ and $v_2$ are two distinct vertices of $G_1$
and $u$ is a vertex of $G_2$. We say that $G^*$ is obtained from $G$ by {\it relocating $G_2$ from $v_2$ to $v_1$}; see Fig. 2.1.

\begin{center}
\includegraphics[scale=.6]{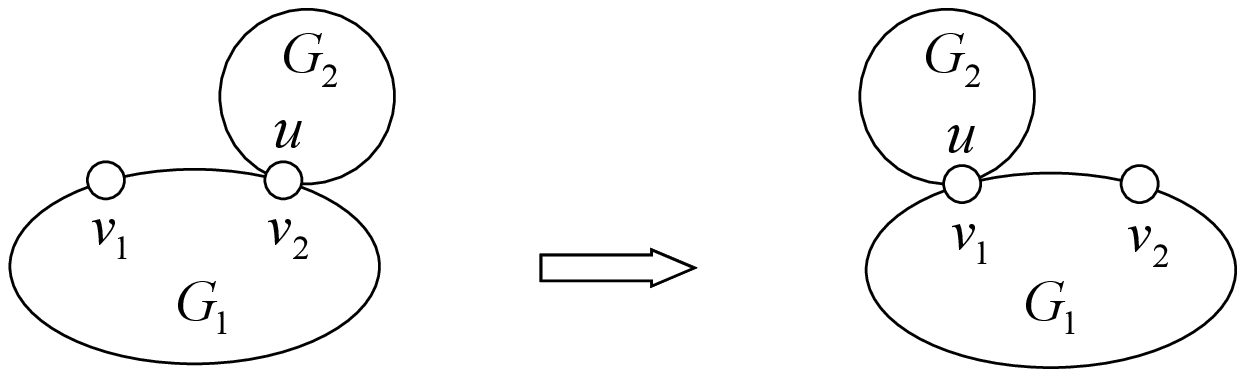}

{\small Fig. 2.1. Relocating $G_2$ from $v_2$ to $v_1$.}
\end{center}

\begin{lemma}{\em \cite{fant}}\label{relocate}
Let $G_1$ be a connected graph containing at least two vertices $v_1$, $v_2$, and let $G_2$ be a nontrivial connected graph containing a
vertex $u$. Let $G=G_1(v_2)\diamond G_2(u)$ and $G^*=G_1(v_1)\diamond G_2(u)$. If there exist a Fiedler vector $x$ of $G$ such that
$x(v_1)\geq x(v_2) \ge 0$ and all vertices in $G_2$ are nonnegatively valuated by $x$, then
$$\alpha(G^*) \leq \alpha(G),$$
with equality if and only if $x(v_1)= x(v_2)=0$, $\sum_{w\in N_{G_2}(u)}x(w)=0$, and $x$ is also a Fiedler vector of $G^*$.
\end{lemma}

Now we investigate the property of the algebraic connectivity of the graph $T(k,l,d)$ listed in Section 1.
Denote $T_d:=T(\lceil \frac{n-d}{2} \rceil, \lfloor \frac{n-d}{2} \rfloor, d)$.

\begin{lemma}{\em \cite{fallatk}} \label{treedia}
Among all trees of order $n$ and diameter $d+1$, the tree $T_d$ is the unique
graph with minimum algebraic connectivity.
\end{lemma}

\begin{lemma}{\em \cite{fant}}\label{ineq1}
{\em (1)}  If $k \ge 2$, $\alpha(T(k,l,d)) > \alpha(T(k-1,l,d+1))$;

\noindent {\em (2)} if $l \ge 2$  $\alpha(T(k,l,d)) > \alpha(T(k,l+1,d+1))$.
\end{lemma}

\begin{lemma}\label{ineq2}
Let $\beta(T_d))=\beta \ge 2$, where $n \ge d+2$ and $n \ge 2\beta+1$.
Then $ 2(\beta-1) \le d \le 2\beta-1$.
Furthermore,
$$\alpha(T_d)\ge \alpha(T_{2\beta-1}),$$
with equality if and only if $d=2\beta-1$.
\end{lemma}

{\it Proof:}
Since $\beta(T_d)=\lfloor \frac{d+2}{2} \rfloor$, then  $\beta \leq \frac{d+2}{2} < \beta+1$,
that is $ 2(\beta-1) \le d <2\beta$.
Hence $ 2\beta-2 \le d \le 2\beta-1$.
It suffices to show that
$\alpha(T_{2\beta-2})\ge \alpha(T_{2\beta-1}).$
%
%Let $T_{2\beta-2}:=T(\Big\lceil \frac{n-(2\beta-2)}{2}\Big\rceil,\Big\lfloor\frac{n-(2\beta-2)}{2}\Big\rfloor,2\beta-2)$.
By Lemma \ref{ineq1} and Lemma \ref{treedia},
we have
\begin{align*}
\alpha(T_{2\beta-2}) &= \alpha \left(T\left(\left\lceil \frac{n-(2\beta-2)}{2}\right\rceil,\left\lfloor\frac{n-(2\beta-2)}{2}\right\rfloor,2\beta-2\right)\right) \\
&>  \alpha\left(T\left(\left\lceil \frac{n-(2\beta-2)}{2}\right\rceil-1,\left\lfloor\frac{n-(2\beta-2)}{2}\right\rfloor,2\beta-1\right)\right)\\
&\geq  \alpha(T_{2\beta-1}).
\end{align*}
\hfill $\blacksquare$

\begin{lemma}\label{ineq3} Let $T_{2\beta_1-1},T_{2\beta_2-1}$ be two trees of
order $n$ with matching number $\beta_1,\beta_2$. If $\beta_1< \beta_2$ and
 $n \ge 2\beta_2+1$,
then
$$\alpha(T_{2\beta_1-1}) >  \alpha(T_{2\beta_2-1}).$$
\end{lemma}

{\it Proof:}
Let $T_{2\beta_1-1}:=T(k,l,2\beta_1-1)$.
By Lemma \ref{ineq1}, we have
\begin{align*}
\alpha[T(k,l,2\beta_1-1)]   & > \alpha[T(k-1,l,2\beta_1)]\\
                          & > \alpha[T(k-1,l-1,2\beta_1+1)]\\
                          & = \alpha[T(k-1,l-1,2(\beta_1+1)-1].
\end{align*}
That is  $\alpha(T_{2\beta_1-1}) >  \alpha(T_{2(\beta_1+1)-1})$.
The result follows by induction on the matching number.
\hfill $\blacksquare$

\begin{lemma} \label{spantree}
Let $G\in\mathcal{M}_{n,\beta}$. Then $G$ contains a spanning tree also with matching number $\beta$ .
\end{lemma}

{\it Proof:}
Let $G\in\mathcal{M}_{n,\beta}$ and let $M$ be a maximum matching of $G$.
Denote $e(G)=|E(G)|$. Clearly, $e(G) \geq n-1$ as $G$ is connected.
The result is certainly true if $e(G)= n-1$, in which case $G$ is a tree.
So we assume that $e(G)> n-1$.

Delete an edge $e_1$ of some cycle of $G$, where $e_1\notin M$, producing a graph $G_1$
such that $\beta(G_1)=\beta(G)$. If $e(G_1)= n-1$, then $G_1$ is a spanning tree of $G$.
If $e(G_1)> n-1$, delete an edge $e_2$ of some cycle of $G_1$, where $e_2\notin M$,
producing a graph $G_2$ such that $\beta(G_2)=\beta(G_1)=\beta(G)$.
We continue the above process until we  arrive at an spanning
tree $G_k$ of $G$ such that $\beta(G_k)=\beta(G_{k-1})=\cdots=\beta(G_1)=\beta(G)$, where $k=e(G)-n+1$.
 \hfill $\blacksquare$

\begin{coro}
Let $G \in \mathcal{M}_{n,\beta}$. If $G$ contains cycles, then $G$ contains a spanning unicyclic graphs with matching number $\beta$.
\end{coro}

{\it Proof:} By Lemma \ref{spantree},
$G$ contains a spanning tree $T$ with matching number $\beta$. The result follows by
adding an edge $e\in E(G)\setminus E(T)$ to the tree $T$.
 \hfill $\blacksquare$

\section{Main results}
We first restrict our discussion to trees with minimum algebraic connectivity.

\begin{theorem}\label{maintree}
Among all trees of order $n$ with matching number $\beta$, where $n \ge 2\beta+1$,
 the tree $T_{2\beta-1}$ is the unique
graph with minimum algebraic connectivity.
\end{theorem}

{\it Proof:}
Here we adopt a similar technique used in the paper \cite{fant}.
If $n=2\beta+1$, the result follows obviously since  $T_{2\beta-1}=P_n$ is the unique minimizing graph
among all connected graphs of order $n$.

Now assume that $n \ge 2\beta+2$.
Let $T$ be a minimizing tree of order $n$ with matching number $\beta$.
If $T$ has exactly two {\it pendant stars} (i.e. the star with maximum possible size centered at a quasi-pendant vertex), then $T=T(k,l,d)$ for some $k,l,d$, where $d \ge 2$.
By Lemma \ref{treedia}, $k=\lceil \frac{n-d}{2} \rceil$, $l=\lfloor \frac{n-d}{2} \rfloor$ and $T=T_d$; by Lemma \ref{ineq2}, $d=2\beta-1$.
The result follows.

Next suppose that $T:=T_0$ has more than two pendant stars, which has $p_0$ pendent vertices
and $q_0$ quasi-pendent vertices. Let $x$ be a Fiedler vector of $T_0$.
First assume $T_0$ is of Type I. Let $N_0=\{v\in V(T_0): x(v)=0\}$.
If $|N_0|\geq 2,$ then there exist at least one zero pendant
star $S$ attached at some vertex say $u$,  and at least one positive
quasi-pendant vertex $w$. Relocating the zero star $S$ at $u$ to $w$,
we will arrive at a new tree $T_1$ such that $\alpha(T_1) < \alpha(T_0) $ by Lemma \ref{relocate}.
Note that $\beta(T_1) \leq \beta(T_0)$ (in fact, $\beta(T_1) < \beta(T_0)$);
otherwise we will get a contradiction to the fact that $T_0$ is minimizing.
If $|N_0|=1$, there exist at least two pendant stars
$S_1,S_2$ both being positive or negative valuated by $x$, attached at $u_1,u_2$ respectively.
Without loss of generality, assume $S_1,S_2$ are both positive and $x(u_1) \ge x(u_2) >0$.
Relocating $S_2$ from $u_2$ to $u_1$, we arrive at a new tree $T_1$ such that
$\alpha(T_1) < \alpha(T_0) $ by Lemma \ref{relocate} and $\beta(T_1) < \beta(T_0)$.

If $T$ is of Type II, then there exist at least two pendant stars
$S_1,S_2$ both being positive or negative valuated by $x$, attached at $u_1,u_2$ respectively.
By the similar way with the case $|N_0|=1$ above,
we also arrive at a new tree $T_1$ such that
$\alpha(T_1) < \alpha(T_0) $ and $\beta(T_1) < \beta(T_0)$.

Repeat the above procession on $T_1$ if $T_1$ has more than two pendant stars and continue
a similar discussion to the resulting tree.
Note that from the $k$-th step to the $(k+1)$-th step, either $p_{k+1}=p_k$ and
$q_{k+1}=q_k-1$, or $p_{k+1}=p_k+1$ and $q_{k+1}=q_k$.
So the above procession will be terminated at the $n$-th step in which the tree $T_n$
has exactly two pendant stars, i.e. $T_n=T(k,l,d)$ for some $k,l,d$, where $d \ge 2$.
Hence
$$\alpha(T)=\alpha(T_0)>\alpha(T_1)> \cdots >\alpha(T_n), ~~\beta(T)=\beta(T_0)>\beta(T_1)>\cdots >\beta(T_n).$$
Therefore, noting that $T_{2\beta-1}$ has matching number $\beta$, by Lemma \ref{treedia} and Lemma \ref{ineq2},
$$\alpha(T_{2\beta-1}) \ge \alpha(T) > \alpha(T_n) \ge \alpha(T_d) \ge \alpha(T_{2\beta(T_n)-1}).$$
However, since $\beta(T_n)<\beta(T)=\beta$, by Lemma \ref{ineq3}, we have
$\alpha(T_{2\beta-1}) <  \alpha(T_{2\beta(T_n)-1})$, a contradiction.
So this case cannot happen and the result follows.\hfill $\blacksquare$

%
%\begin{lemma} {\em \cite{fied1}}\label{spansubgraph}
%Let $H$ be a spanning subgraph of $G$, Then
%$$\alpha(H)\leq \alpha(G).$$
%\end{lemma}

\begin{theorem}\label{main}
Let $G\in\mathcal{M}_{n,\beta}$. Then $G$ is minimizing in $\mathcal{M}_{n,\beta}$
if and only if $G= T_{2\beta-1}$.
\end{theorem}

{\it Proof:}
If $\beta=1$, the result holds as $T_1=S_{1,n-1}$ is the unique graph of matching
number $1$ for $n\geq 2$ and $n\neq 3$.
When $n=3$, there are exactly two graphs: $S_{1,2}$ and the triangle $C_3$, both having matching number $1$.
Since $\alpha(S_{1,2})<\alpha(C_3)$, the result also holds in this case.

Assume that $\beta\geq 2$. If $n=2\beta$, the result surely holds as $P_{2\beta}$ is the unique
minimizing graph. So suppose that $n\geq 2\beta+1$ in the following.
Let $G$ be a minimizing graph in $\mathcal{M}_{n,\beta}$.
Then $G$ contains a spanning tree $T$ with matching number $\beta$ by Lemma \ref{spantree}.
Furthermore, by Theorem \ref{maintree},
$$\alpha(G) \ge \alpha(T) \ge \alpha(T_{2\beta-1}).\eqno(3.1)$$
 Hence $\alpha(G) = \alpha(T) = \alpha(T_{2\beta-1})$, which implies that  $T= T_{2\beta-1}$ also by Theorem \ref{maintree}.

 We claim that $G=T_{2\beta-1}$; otherwise $E(G) \backslash E(T_{2\beta-1}) \ne \emptyset$.
 Let $x$ be a unit Fiedler vector of $G$.
Then
\begin{align*}
\alpha(G) &= \sum_{uv \in E(G)}[x(u)-x(v)]^2\\
& =\sum_{uv \in E(T_{2\beta-1})}[x(u)-x(v)]^2+\sum_{uv \in E(G)\backslash E(T_{2\beta-1})}[x(u)-x(v)]^2\\
& \geq \sum_{uv \in E(T_{2\beta-1})}[x(u)-x(v)]^2 \\
& \geq \alpha(T_{2\beta-1}).
\end{align*}
Since $\alpha(G) =\alpha(T_{2\beta-1})$, then $x$ is also  a Fiedler vector of $T_{2\beta-1}$, and $x(u)=x(v)$ for each edge
$uv \in E(G)\backslash E(T_{2\beta-1})$.
 By Lemma \ref{fiedler2}, whenever $T_{2\beta-1}$ is of Type I or Type II,
 $u,v$ should be both the pendent vertices lying in a same pendent star.
However, in this case  $\beta(T_{2\beta-1}+uv) > \beta=\beta(G)$ for any $uv \in E(G)\backslash E(T_{2\beta-1})$; a contradiction.

The sufficiency results follows from the discussion of (3.1).
\hfill $\blacksquare$
%
%\begin{lemma}{\em \cite{Gallai}} \label{cover}
%For every graph $G$ of order $n$ containing no isolated vertices,
%$$\be(G) + \nu(G) = n.$$
%\end{lemma}

As a byproduct, we get the following result on edge covering number.

\begin{coro}\label{main2}
Let $G\in\mathcal{C}_{n,\gamma}$. Then $G$ is minimizing in $\mathcal{C}_{n,\gamma}$
if and only if $G= T_{2(n-\gamma)-1}$.
\end{coro}

%\begin{theorem}\label{main2}
%Let $G$ be a connected graph of order $n$ and edge covering number $\nu$, where $n \leq 2\nu$.
%Then $G$ is minimizing graph
%if and only if $G\cong K_{1,n}$, for $\nu=n-1$, or  $G\cong T_{2(n-\nu)-1}$, for $\frac{n}{2}\leq\nu\leq n- 2$.
%\end{theorem}

{\it Proof:}
The result follows by Theorem \ref{main} and the fact that $\beta(G)+\gamma(G)=n$.
 \hfill $\blacksquare$

\begin{lemma}{\em \cite{kirkn}} \label{bound}
Suppose that $d \ge 3$, $ k \ge 1$, $l \ge 1$ and $n:=k+l+d-1$.
Then
$$\alpha(T(k,l,d-1)) \ge \left(\frac{nd}{4}-\frac{2n+d^2-4d-5}{8}\right)^{-1}.$$
\end{lemma}

\begin{coro}
Let $G\in \mathcal{M}_{n,\beta}$.
Then $$\alpha(G) \ge \frac{8}{-4\beta^2+4\beta(n+2)-2n+5}.$$
\end{coro}

{\it Proof:}
By Theorem \ref{main}, $\alpha(G) \ge \alpha(T_{2\beta-1})$.
If $\beta=1$, surely $$\alpha(T_{2\beta-1})=\alpha(T_1)=1 > \frac{8}{-4\beta^2+4\beta(n+2)-2n+5}=\frac{8}{2n+9}.$$

If $\beta \ge 2$ and $n =2\beta$, noting that in this case $T_{2\beta-1}= P_n$,
$$\alpha(P_{n)}=2\left(1-\cos\frac{\pi}{n}\right)> \frac{8}{-4\beta^2+4\beta(n+2)-2n+5}=\frac{8}{n^2+2n+5}.$$

If $\beta \ge 2$ and $n \geq2\beta+1$,
the result follows by taking $d=2\beta$ in Lemma \ref{bound}.
 \hfill $\blacksquare$

Similarly we have the following corollary.

\begin{coro}
Let $G\in \mathcal{C}_{n,\gamma}$.
Then $$\alpha(G) \ge \frac{8}{-4\gamma^2+4\gamma(n-2)+6n+5}.$$
\end{coro}

\end{document}